\documentclass[11pt]{amsart}

\usepackage{esint}
\usepackage{color}
\usepackage{hyperref}
\hypersetup{bookmarksdepth=3}

\newtheorem{thm}{Theorem}[section]
\newtheorem{lm}[thm]{Lemma}

\theoremstyle{definition}
\newtheorem{df}[thm]{Definition}
\newtheorem*{df*}{Definition}

\theoremstyle{remark}
\newtheorem{rem}[thm]{Remark}
\newtheorem*{rem*}{Remark}

\numberwithin{equation}{section}

\newcommand{\ci}[1]{_{ {}_{\scriptstyle #1}}}

\newcommand{\cB}{\mathcal{B}}
\newcommand{\cD}{\mathcal{D}}

\newcommand{\cX}{\mathcal{X}}
\newcommand{\cL}{\mathcal{L}}
\newcommand{\cM}{\mathcal{M}}

\newcommand{\f}{\varphi}

\newcommand{\e}{\varepsilon}

\newcommand{\R}{\mathbb{R}}
\newcommand{\Z}{\mathbb{Z}}

\newcommand{\E}{\mathbb{E}}

\newcommand{\PP}{\mathbb{P}}

\newcommand{\om}{\omega}

\newcommand{\ch}{\operatorname{ch}}
\newcommand{\p}{\partial}

\newcommand{\mfA}{\mathfrak{A}}

\newcommand{\1}{\mathbf{1}}

\newcommand{\wt}{\widetilde}

\newcommand{\cz}{Calder\'{o}n--Zygmund\ }

\newcommand{\spn}{\operatorname{span}}

\newcommand{\La}{\langle }
\newcommand{\Ra}{\rangle }

\newcommand{\ssd}{{\scriptscriptstyle\Delta}}
\newcommand{\sd}{{\scriptstyle\Delta}}
\newcommand{\bu}{\mathbf{u}}
\newcommand{\bn}{\mathbf{n}}

\newcommand{\bw}{\mathbf{w}}
\newcommand{\bff}{\mathbf{f}}

\newcommand{\be}{\mathbf{e}}

\newcommand{\av}[2]{\langle #1\rangle_{_{\scriptstyle #2}}}

\newcommand{\bT}{\mathbf{T}}

\newcommand{\fdot}{\,\cdot\,}

\def\cyr{\fontencoding{OT2}\fontfamily{wncyr}\selectfont}
\DeclareTextFontCommand{\textcyr}{\cyr}
\newcommand{\sha}[0]{\ensuremath{\mathbb{S}
}}

\newenvironment{entry}
{\begin{list}{X}%
  {%
      \setlength{\labelwidth}{55pt}%
      \setlength{\leftmargin}{\labelwidth}
      \addtolength{\leftmargin}{\labelsep}%
   }%
}%
{\end{list}}

\renewcommand{\labelenumi}{(\roman{enumi})}

\newcounter{vremennyj}

\newcommand\cond[1]{\setcounter{vremennyj}{\theenumi}\setcounter{enumi}{#1}\labelenumi\setcounter{enumi}{\thevremennyj}}

\begin{document}

\title
{A Bellman function proof of the $L^2$ bump conjecture
}
\author[F.~Nazarov]{Fedor Nazarov}
\author[A.~Reznikov]{Alexander Reznikov}
\author[S.~Treil]{Sergei Treil}
\author[A.~Volberg]{Alexander Volberg}
\thanks{Partially supported by NSF }


\makeatletter
\@namedef{subjclassname@2010}{
  \textup{2010} Mathematics Subject Classification}
\makeatother




%
%


\maketitle
\begin{abstract}
We approach the problem of finding the sharp sufficient condition of the boundedness of all two weight \cz operators. We solve this problem in $L^2$ by writing a formula for a Bellman function of the problem.
\end{abstract}

\section{Introduction}
\label{intro}

\subsection{Preliminaries}
In this paper we give a simple Bellman function solution of the so-called ``bump conjecture'' for the two weight estimates of the singular integral operators.

The original (still open) question about two weight estimates for the singular integral operators is to find a necessary and sufficient condition on the weights $u$ and $v$ such that a \cz operator $T: L^p(u)\to L^p (v)$ is bounded, i.e.~the inequality
\begin{align}
\label{two-weight_01}
\int |Tf|^p v dx \le C \int |f|^p u dx \qquad \forall f \in L^p(u)
\end{align}
holds.

In the one weight case $v=u$ the famous Muckenhoupt condition is necessary and sufficient for \eqref{two-weight_01}
\begin{align}
\tag{$A_2$}
\sup_I \left(|I|^{-1}\int_{I} v dx \right) \left(|I|^{-1}\int_{I} v^{-p'/p} dx \right)^{p/p'} <\infty
\end{align}
where the supremum is taken over all cubes $I$. More precisely, this condition is sufficient for all \cz operators, and is also necessary for classical (interesting) \cz operators, such as Hilbert transform, Riesz transform (vector-valued, when all Riesz transforms are considered together), Beurling--Ahlfors operator.

The inequality \eqref{two-weight_01} is equivalent to the boundedness of the  operator $M_{v^{1/p}} T M_{u^{-1/p}}$ in the non-weighted $L^p$; here $M_\f$ is the multiplication operator, $M_\f f = \f f$. Denoting $w=u^{-p'/p}$ we can rewrite the problem in the symmetric form as the $L^p$ boundedness of $M_{v^{1/p} }T M_{w^{1/p'}}$.

So the problem can be stated as: \emph{Describe all weights (i.e.~non-negative  functions) $v$, $w$ such that the operator $M_{v^{1/p} }T M_{w^{1/p'}}$ is bounded in (the non-weighted) $L^p$}

Note, that this symmetric formulation is  more general than \eqref{two-weight_01}, because in \eqref{two-weight_01} it is usually assumed that $u$ and $v$ are locally integrable, but  for \eqref{two-weight_01} to hold for interesting operators (Hilbert Transform, vector Riesz Transform, Beurling--Ahlfors Transform, etc.) the function $1/u$ also has to be locally integrable.

For the interesting operators the following two weight analogue of the $A_p$ condition is necessary for the boundedness of the operator $M_{v^{1/p} }T M_{w^{1/p'}}$:
\begin{align}
\sup_I \left(|I|^{-1}\int_{I} v dx \right) \left(|I|^{-1}\int_{I} w dx \right)^{p/p'} <\infty
\end{align}
or in the symmetric form
\begin{align}
\label{2weightA_p}
\sup_I \left(|I|^{-1}\int_{I} v dx \right)^{1/p} \left(|I|^{-1}\int_{I} w dx \right)^{1/p'} <\infty
\end{align}

Simple counterexamples show that this condition is not sufficient for the boundedness. So a natural way to get a sufficient condition is to replace the $L^1$ norms of $v$ and $w$ in \eqref{2weightA_p}  (or the $L^p$ and $L^{p'}$ norms of $v^{1/p}$ and $w^{1/p'}$) by some stronger Orlicz norms (``bumping'' the $L^p$ norms).

Namely, given a Young function $\Phi$ and a cube $I$ one can consider the normalized on $I$ Orlicz space $L^\Phi(I)$ with the norm given by
\begin{align*}
\|f\|\ci{L^\Phi(I)} := \inf\left\{ \lambda>0 : \int_I \Phi\left( \frac{ f(x)}{\lambda } \right) \frac{dx}{|I|} \le 1 \right\}.
\end{align*}
And it was conjectured (for $p=2$) that  if the Young functions $\Phi_1$ and $\Phi_2$ are integrable near infinity,
\begin{equation}
\label{inte}
\int^\infty \frac{dx}{\Phi_{i}(x)} <\infty,\,\, i=1,2\,,
\end{equation}
then the condition
\begin{align}
\label{bump_01}
\sup_I \| v\|\ci{L^{\Phi_1}(I)} \| w\|\ci{L^{\Phi_2}(I)} < \infty
\end{align}
implies that  for any bounded \cz operator $T$
the operator $M_{v^{1/2}} T M_{w^{1/2}}$ is bounded in $L^2$. Usually in the literature a more complicated (although equivalent) form of this conjecture was presented, but at least in the case $p=2$ condition \eqref{bump_01} seems more transparent.%
\footnote{The bump condition was also  stated for $p\ne 2$, but  in this paper we only deal with the case $p=2$.}

Condition \eqref{bump_01} was considered in numerous papers in the attempt to prove its universal sufficiency for {\bf all} \cz operators. The reader can find beautiful approaches in  \cite{CU-Ma-Pe}, \cite{CU-Ma-Pe07}, \cite{CU-Pe99},  \cite{CU-Pe00},  \cite{CU-Pe-pisa},  \cite{Le},  \cite{Pe94JL}, \cite{PeMax},  where partial results for some \cz operators were proved (note that \cite{PeMax} is about maximal operator and not about \cz operators). Finally in \cite{Le1} the sufficiency of bump condition for all \cz operators to be bounded was fully proved (and even generalized to all $p\in (1, \infty)$), although in formally less general situation of the weights.

Simultaneously and independently the first version of the present paper \cite{NRTV1}
appeared. Slightly earlier the sketch of the approach (with the main  ideas but without much details) was circulated as \cite{NRTV}. A diligent reader will recognize that the approaches in the present paper (and thus in earlier versions \cite{NRTV}, \cite{NRTV1}) and in Lerner's paper \cite{Le1} are very different, but still have something in common. This  very important  common point is the ``coupling by the same cube" feature. It is the main winning idea of \cite{Le1}. And it is the feature of the present paper (and \cite{NRTV}, \cite{NRTV1}) as Bellman function approach automatically should have this feature.

\subsection{What is done in the paper}
Formally, in this paper we prove the $L^2$ case of the bump conjecture using Bellman function method. As it is now well-known, a general \cz operator can be represented as an average of dyadic  shift and paraproducts, so it is sufficient to prove the estimates for such operators, and that is exactly what is done in the paper.

However we think that the results obtained in the paper that were used to prove the bump conjecture are of significant  interest by themselves; probably they are even more interesting than the solution of the bump conjecture.


Let us shortly describe what is done in the paper.

\begin{itemize}
\item First, the Orlicz norm is not easy to work with.  We introduce a lower bound for the Orlicz norm,  which gives a more tractable, in our opinion,  way to ``bump'' the averages.  In particular, it allows us to apply the Bellman function method.
\item The application of the Bellman function method is by now standard.  The novelty of the argument belongs to the fact that Bellman function now is defined on an infinite-dimensional space.
\item The estimates for the Haar shifts and for the paraproducts are reduced to two embedding theorems, so   the operators are \emph{constructively} factorized through $\ell^2$. This essentially means that the bump condition is a rather rough one, since in  more delicate two weight situations no such factorization appeared to be possible (at least no factorization was found),  see, for example,  \cite{NTV99}, \cite{CUReVo}.
\item Namely, it is known, see \cite{N83}, that in the case of  power bumps ($\Phi(t) =t^{1+\e}$) one can insert a Muckenhoupt $A_p$ weight between $w^{-1}$ and $C v$, so the boundedness follows immediately. For finer bumps such insertion of $A_p$ weights is impossible, but the constructive factorization through $\ell^2$ can be considered the next  best thing.
\item Finally, the main estimates can be directly extended to general (non-homogeneous) martingale settings, and can be used in more general situations. In particular,  a word by word extension of our results gives the proof of the bump conjecture for the \cz operators on geometrically doubling metric spaces (equipped with a doubling measure).
Indeed, since random ``dyadic'' lattices can be constructed on geometrically doubling metric spaces, and representation of \cz operators on such spaces as an average of Haar shifts is now known (see, for example, \cite{NRV}), everything follows from our results (see Theorem \ref{t:shift-para} below).
\end{itemize}

\noindent{\bf Acknowledgement.} We are grateful to the referee for several helpful remarks.

\subsection{Orlicz norms and distribution functions}
\label{orlicz}

Orlicz norm is not very convenient to work with, so we would like to replace it by something more tractable.

\subsubsection{A lower bound for the Orlicz norm} We start with the remark that notation $\int_0 f(t)dt <\infty$ means that the function is integrable near zero. Similarly, $\int^\infty f(t)dt<\infty$ means that the function is integrable near infinity.

Let $\Phi$ be a continuous non-negative increasing convex function such that $\Phi(0)=0$ and
$\int^{+\infty}\frac{dt}{\Phi(t)}<+\infty$. Define $\Psi(s)$ parametrically by $\Psi(s)=\Phi'(t)$ when
$s=\frac{1}{\Phi(t)\Phi'(t)}$ ($t>0$). Then $\Psi(s)$ is positive and decreasing for
$s>0$ and $s\Psi(s)$ is increasing. Moreover $\int_0\frac{ds}{s\Psi(s)}<+\infty$. Indeed, using our parametrization we can rewrite the last integral as
$$
\int^{+\infty}\left(\frac{1}{\Phi(t)}+\frac{\Phi''(t)}{\Phi'(t)^2}\right)\,dt\,.
$$
The first integral converges by our assumption and the second integrand has a bounded near $+\infty$
antiderivative $\frac{-1}{\Phi'(t)}$.

Let $w\ge 0$ on  $I\subset \R^n$. Define the normalized distribution function $N$ of $w$ by
\begin{align}
\label{distr_fn_01}
N(t)=N_I^w(t)=\frac{1}{|I|} \left| \{x\in I:w(x)>t\} \right|
\end{align}

\begin{lm}
\label{l:Orl_LB}
Let $\Psi : (0,1]\to \R_+$ be a decreasing function such that the function $s\mapsto s\Psi(s)$ is increasing. Let $\Phi$ be a Young function and let
\begin{align*}
\Psi(s) \le C \Phi'(t) \qquad \text{where} \quad s = \frac{1}{\Phi(t) \Phi'(t)}
\end{align*}
for all sufficiently large $t$. Then for $N=N_I^w$
\begin{align}
\label{n(N)}
\bn\ci\Psi (N) := \int_0^\infty N(t)\Psi(N(t))\,dt\le C\|w\|\ci{L^\Phi(I)}\,.
\end{align}
\end{lm}

\begin{proof}
The left hand side scales like a norm under multiplication by constants, so it is enough to
show that if $\|w\|_{L^\Phi(I)}\le 1$, i.e.,
\[
\frac{1}{|I|}\int_I\Phi(w)=\int_0^\infty N(t)\Phi'(t) \,dt \le 1
\]
then $\bn\ci\Psi(N)$ is bounded by a constant. Since $s\Psi(s)$ increases, we may have
trouble only at $+\infty$ It is clear that it suffices to estimate the
integral over the set
where $\Psi(N(t))>\Phi'(t)$ but since $\Psi$ is decreasing this means that $N(t)\le C/{ (\Phi(t)\Phi'(t)) }$, so we get
at most $\int^{+\infty}\Phi(t)^{-1} dt$ and we are done.
\end{proof}

\begin{rem*}
In fact, for sufficiently regular $\Phi$, the converse inequality $\|w\|\ci{L^\Phi}\le C \int_0^\infty N(t)\Psi(N(t))\,dt$ holds for any positive decreasing integrable $N$.  To see this, let us consider the family of $\Phi$'s such that $\Phi(t) = t\rho(t)$ and $\rho$ is  monotonically increasing and ``logarithmically  concave" in the sense that $\frac{t\rho'(t)}{\rho(t)}$ decreases monotonically when $t\rightarrow\infty$. We also assume of course that $\lim_{t\rightarrow\infty} \rho(t)=\infty$ and that $\rho(t)\ge 1$.
Let
$G(t):=N(t) \Psi(N(t))$. When t goes to infinity, $N$ is monotonically decreasing to zero, and hence $G$ is also monotonically decreasing (as $s\Psi(s)$ increases  near zero).

 Put $s^{-1} =  \Phi(t)\Phi'(t)\asymp t\rho^2(t)$ (just because $\Phi'(t)\asymp \rho(t)$ by our ``logarithmic concavity" of $\rho$ assumption). Hence $s\ge c_1 (t\rho^2(t))^{-1}$. Now  $\Psi$ is decreasing by definition, and this implies

 \begin{equation}
 \label{first}
 \Psi(c_1( t\rho^2(t))^{-1}) \ge \Phi'(t)\asymp \rho(t)\ge c_2 \rho(t)\,.
 \end{equation}

 We now ask an addition to ``logarithmic concavity", namely:
 \begin{equation}
 \label{rho}
 \frac{t\rho'(t)}{\rho(t)}\log\rho(t)\rightarrow 0\,.
 \end{equation}
 Denote $r(x)=\log(\rho(e^x))$.  We required  at the beginning that $\lim_{x \rightarrow\infty}r(x)=\infty$. The last inequality says in particular that $r'(x)=o(1) r(x)^{-1}$, and therefore, $r'$ tends to zero at infinity. Thus $r(x)\le \frac{x}{3}$ for all large $x$. Keeping this in mind we continue.

Set
$
u=\frac{t\rho^2(t)}{c_1}.
$
Then
$$
t=\frac{c_1u}{\rho^2(t)}.
$$
Thus, since $\rho$ is an increasing to infinity function and we assume that $t$ is sufficiently big, we get $\rho^2(t)\ge c_1$. Therefore,
$$
t\le u,
$$
and, thus,
$$
t=\frac{c_1u}{\rho^2(t)} \ge c_1 \frac{u}{\rho^2(u)}.
$$
Hence, using \eqref{first}, we get
\begin{equation}
\label{second}
\Psi(u^{-1}) \ge c_2 \rho(t) \ge c_2\rho(c_1\frac{u}{\rho^2(u)})\,. 
\end{equation}

Next, we will prove the following inequality. Recall that $r(x)=\log(\rho(e^x))$. We claim that
$$
\Delta(x)=r(x)-r(x-2r(x)-c_0) \le C.
$$
In fact, by the mean value theorem we have for certain $\xi\in (x-2r(x)-c_0, x)$
$$
\Delta(x) = (2r(x)+c_0) r'(\xi) = (2r(x)+c_0) \frac{\rho'(e^\xi)}{\rho(e^\xi)}e^\xi.
$$

Since we assumed that $t \frac{\rho'(t)}{\rho(t)}$ is monotonically decreasing, we get (now using \eqref{rho} in the second comparison below):
\begin{multline*}
\Delta(x) \le (2r(x)+c_0) \frac{\rho'(e^{x-2r(x)-c_0})}{\rho(e^{x-2r(x)-c_0})} e^{\rho(e^{x-2r(x)-c_0})} = \\ =
(2r(x)+c_0) \frac{o(1)}{\log \rho(e^{x-2r(x)-c_0})} =(2r(x)+c_0) \frac{o(1)}{r(x-2r(x)-c_0)} = \\=
 o(1) \left(\frac{2r(x)+c_0}{r(x-2r(x)-c_0)} -2 + 2\right) = o(1) \frac{2\Delta(x)+c_0}{r(x-2r(x)-c_0)} + o(1).
\end{multline*}
Finally, we use that $r(x-2r(x)-c_0)$ is separated from zero when $x$ is big. Thus
$$
\Delta(x)\le \Delta(x)o(1) + o(1).
$$
This immediately implies $\Delta(x)=o(1)$ when $x\to \infty$, and thus
\begin{equation}
\label{third}
\Delta(x)\le C\,.
\end{equation}

Let us now write what does it mean. In fact, by the definition of $r$ and by \eqref{third}, we can conclude that
$$
C\ge r(x)-r(x-2r(x)-c_0) = \log\frac{\rho(e^x)}{\rho(e^{x-2\log\rho(e^x)-c_0})} = \log\frac{\rho(e^x)}{\rho(\frac{e^x}{c_3\rho^2(e^x)})}.
$$
Thus, we get  for all large $u$:
$$
\rho(u)\le c_4 \rho(\frac{u}{c_3 \rho^2(u)})\,.
$$
We chose $c_3 = c_1^{-1}$ and plug the above inequality into \eqref{second}. Then we finally get
$$
\Psi(u^{-1}) \ge c_5\rho(u)
$$

If $N\Psi(N)= G$ then $c_6  G \ge N\rho(\frac1N)$ by the previous inequality. Therefore, $N\le c_6G$ (we assumed that $\rho\ge 1$), and  $N\le\frac{c_6G}{\rho(\frac1N)} \le \frac{c_6G}{\rho(\frac1{c_6G})}$.  And we can continue the previous estimate: $N(t) \le  \frac{c_6G(t)}{\rho(\frac1{c_6G(t)})}\le \frac{c_6G(t)}{\rho(t)}$.  We used here the fact that the integrability and monotonicity of $G$ implies that $G(t)= o(\frac1t)$, in particular, $G(t)<\frac1{c_6t}$ for large $t$. But we already mentioned that $\Phi'(t) \le c_7 \rho(t)$. Combining the last two inequalities, we get
$N(t)\Phi'(t) \le c_6c_7 G(t)$, and we just obtained that $\int_0^{\infty} N(t) \Phi'(t) dt \le c_4 c_5$.
\end{rem*}

\subsubsection{Examples} In the above section only the behavior of $\Phi$ at $+\infty$ (equivalently, the behavior of $\Psi $ near $0$) was important, so we will concentrate our attention there.

Let $\Phi(t) = t (\ln t)^\alpha$, $\alpha>1$ near $\infty$. Then
\[
\Phi'(t) \sim (\ln t)^{\alpha}, \qquad \Phi(t)\Phi'(t) \sim t (\ln t)^{2\alpha} ,
\]
so $\Psi(s) := (\ln (1/s))^\alpha$ satisfies the assumptions of Lemma \ref{l:Orl_LB}: to see that we  notice
\[
\ln(\Phi(t)\Phi'(t) )\sim \ln t.
\]

If $\Phi(t) = t \ln t (\ln\ln t)^\alpha$, $\alpha>1$, then
\[
\Phi'(t)\sim  \ln t (\ln\ln t)^\alpha, \qquad \Phi(t)\Phi'(t) \sim t (\ln t)^2 (\ln\ln t)^{2\alpha}
\]
and $\Psi(s) = \ln(1/s) (\ln\ln (1/s))^\alpha$ works. because again $\ln(\Phi(t)\Phi'(t) )\sim \ln t$.

Note that in both examples $\int_0 (s\Psi(s))^{-1} ds <\infty$.

The examples of Young functions with higher order logarithms are treated similarly.

\subsection{Main result}
Let $\Psi_1, \Psi_2:(0,1] \to \R_+$ be as above, i.e. for $i=1,2$, $\Psi_{i}$ is decreasing, $s \mapsto s\Psi_{i}(s)$ is increasing and
\[
\int_0^1 \frac{ds}{s\Psi_{i}(s)} <\infty.
\]
Recall that for a weight $w$ the normalized distribution function $N_I^w$ is defined by \eqref{distr_fn_01}
\begin{thm}
\label{t:main_01}
Let the weights $v$, $w$ satisfy
\begin{equation}
\label{bump_02}
\sup_I \bn\ci{\Psi_1}(N_I^v) \bn\ci{\Psi_2}(N_I^w) < \infty;
\end{equation}
here the supremum is taken over all cubes $I$, and $\bn\ci{\Psi}$ is defined by \eqref{n(N)}.

Then for any bounded \cz operator $T$ the operator
\[
M_{v^{1/2}} T M_{w^{1/2}}
\]
is bounded in $L^2$.
\end{thm}

\section{Reductions: Haar shifts, paraproducts and embedding theorems}
First, let us reduce the problem to its dyadic (martingale) analogue, i.e.~to the estimates of the so-called Haar shifts and paraproducts.

Since a bounded \cz operator can be represented as a weighted average (over the random dyadic grids) of Haar shifts and paraproducts and their adjoints, where the weights decay exponentially in complexity of the Haar shifts, it is sufficient to get the estimates for the Haar shifts that grow sub-exponentially (for example, polynomially) in the complexity of the shifts and the estimates for the paraproducts (there is no complexity of the paraproducts).

The estimates for each operators will be in turn factored into  two embedding theorems, and these embedding theorems are proved in this paper.

The embedding theorems and so the estimates of the Haar shifts and paraproducts hold in very general martingale settings,

\subsection{General setup}
\label{setup}

Consider a measure space $X$ with $\sigma$-finite measure  $\mu$ let $\cL_k=\{I_j^k\}_{j}$, $k\in\Z$ (or $k\in \Z_+$) be  partitions of $X$ into disjoint sets $I_j^k$, $0<\mu(I_j^k)<\infty$.

We assume that the partition $\cL_{k+1}$ is a refinement of $\cL_k$.

Let $\mfA$ be the $\sigma$-algebra generated by all the partitions $\cL_k$. In what follows all functions on $X$ we consider will be   assumed to be $\mfA$-measurable.

With respect to this $\sigma$-algebras we can define martingale averaging operators $\E_k$, and martingale difference operators $\Delta^n_k:= -\E_k +\E_{k+n}$.

We adapt the following notation.

\begin{entry}
\item[$\ch I$] The collection of children of  $I\in\cL$, i.e.~if $I\in \cL_n$ then $\ch I =\{J\in\cL_{n+1} : J\subset I\}$.
\item[$\ch_k I$] The collection of children of the order $k$ of  $I\in\cL$; $\ch_0(I) =\{I\}$, $\ch_{k+1}(I) =\{\ch(J):J\in\ch_k(I)\}$.
\item[$\La f \Ra \ci I$, $\fint_I f$] The average of $f$ over $I$, $\La f \Ra \ci I = \mu(I)^{-1} \int_I f(x) d\mu(x)$;
\item[$E\ci I$] The averaging operator, $E\ci I f := \La f \Ra\ci I \1\ci I$;
 note that \\ $E_k = \sum_{I\in\cL_k} E\ci I$.
\item[$\Delta\ci I$] Martingale difference operator, $\Delta\ci I:= -E\ci I + \sum_{J\in\ch(I)} E\ci J$; note that  $\Delta_k = \sum_{I\in \cL_k} \Delta\ci I$.
\item[$\Delta^n\ci I$] Martingale difference operator of order $n$,
\[
\Delta^n\ci I:= -E\ci I + \sum_{J\in\ch_n(I)} E\ci J.
\]
\end{entry}

Since the measure $\mu$ is assumed to be fixed we sometimes will be using $|E|$ for $\mu(E)$ and $dx$ for $d\mu(x)$. We also will be using $L^2$ for $L^2(\mu)$

The prototypical  example is $X= \R$ or $\R^d$ with $\cL$ being a dyadic lattice  $\cD$.

\subsection{Haar shifts}
\label{DS2w}

 \begin{df}
\label{df-Haar2}
A Haar shift $\sha$ of complexity $n$ is given by
\[
\sha f=\sum_{I\in\cD} \sha\ci I \Delta\ci I^n f,
\]
where the operators $\sha\ci I$  act on $\Delta^n\ci I L^2$ and can be represented as integral operators with kernels $a\ci I$, $\|a\ci I\|_\infty\le |I|^{-1}$. The latter means that for all $f, g\in \Delta\ci I^n L^2$
\[
\La \sha \ci I f, g\Ra = \int_I\int_I a\ci I(x, y) f (y) g(x) dx dy.
\]
\end{df}

This is a slightly more general definition than the one in  \cite{HPTV}, but only the estimate $\|a\ci I\|_\infty \le |I|^{-1}$ is essential for our construction. Note also that according to the definition in \cite{HPTV} the complexity of the corresponding shift is $n-1$, not $n$, which really does not matter; we just find our definition of complexity a bit more convenient.

The estimate $\|a\ci I\|_\infty \le |I|^{-1}$ means that the operators $\sha\ci I$ are ``$L^1\times L^1$ normalized'', meaning that
\begin{align}
\label{L^1xL^1_normalization}
|\La\sha\ci I  f,  g \Ra | \le |I| \frac{\|f\|_1}{|I|} \frac{\|g\|_1}{|I|} \qquad \forall f, g \in \Delta^n\ci I L^2
\end{align}
Haar shifts of complexity $1$ are simply ``$L^1\times L^1$ normalized'' martingale transforms; martingale transform here means in particular that the subspaces $\Delta\ci I$ are orthogonal, and $\sha$ can be represented as an orthogonal sum of the operators $\sha\ci I$.

A Haar shift of complexity $n\ge2$ is not generally a martingale transform, meaning that the subspaces $\Delta\ci I^n$  generally intersect, so $\sha$ does not split into direct sum of $\sha\ci I$.

However, if one goes with step $n$, then the corresponding operator is a martingale transform, so a Haar shift of complexity $n$ can be represented as a sum of $n$ Haar shifts of complexity $1$. Namely, for $k=1, 2, \ldots, n-1$ define
\[
\cL^k =\{I\colon I\in\cL_{k+nj}, j\in\Z\},
\]
and let
\[
\sha_k = \sum_{I\in\cL^k} \sha\ci I.
\]
Then $\sha = \sum_{k=0}^{n-1} \sha_k $ and each $\sha_k$ is a Haar shift of complexity $1$ with respect to the lattice $\cL^k$.

\begin{rem*}
Therefore, uniform estimate for the Haar shifts of complexity $1$ (i.e.~for the ``$L^1\times L^1$ normalized'' martingale transforms) gives the linear in complexity estimate for the general Haar shifts. Notice that the estimate does not depend on the number of children.
\end{rem*}

\subsection{Paraproducts}

Given the lattice $\cL$ and a locally integrable function $b$, the paraproduct $\Pi = \Pi_b =\Pi_b (\cL)$ is defined as
\[
\Pi f :=  \sum_{I\in\cL} (E\ci I f) (\Delta\ci I b).
\]
The necessary and sufficient condition for the paraproduct to be bounded is that
\[
\sup_{J\in\cL} |J|^{-1} \sum_{I\in \cL: I\subset J} \|\Delta\ci I b\|_2^2 <\infty.
\]
In the case of  dyadic lattice in $\R^d$ or, more generally in the homogeneous situation, when
\[
\inf_{J\in\cL} \ \inf_{I\in \ch(J)} \frac{|I|}{|J|} >0
\]
this condition is equivalent to $b$ belonging to the corresponding martingale BMO space $\text{BMO}\ci\cL$

\subsection{Reduction to the martingale case.}
\label{redu1}

To reduce the problem to the martingale case we  use  the following result that can be found in \cite{H} and \cite{HPTV}:
\begin{thm}
\label{deco}
Let $T$ be a \cz operator in $\R^d$. There exists a probability space $(\Omega, \PP)$ of dyadic lattices $\cD_{\om}$, such that
\[
T =C\left( \int_\Omega \sum_{n=1}^\infty 2^{-\e n} \, \sha_n(\om) d\PP(\om) + \int_\Omega ( \Pi^1 (\om) +  (\Pi^2(\om))^*) d \PP(\om) \right)\,,
\]
where $\sha_n(\om)$ are Haar shifts of complexity $n$ with respect to the lattice $\cD_\om$, $\Pi^{1, 2} (\om)$ are the paraproducts with respect to the lattice $\cD_\om$, $\|\Pi^{1, 2}(\om)\|\le 1$.

The constants $C$ and $\e$ depend  on $d$, $\|T\|$ and \cz parameters of the kernel of $T$.
\end{thm}

Theorem \ref{deco} implies immediately that
 the main theorem (Theorem \ref{t:main_01}) follows from the theorem below.
\begin{thm}
\label{t:shift-para}
Let the weights $v$, $ w$ satisfy the assumptions of Theorem \ref{t:main_01}. Then
\begin{enumerate}
	\item For all Haar shifts $\sha$ of order $1$ the operators $M_{v^{1/2}} \sha M_{w^{1/2}}$ are uniformly bounded in $L^2$, $\|M_{v^{1/2}} \sha M_{w^{1/2}}\| \le C$, where $C$ depends on $\Psi_2$, $\Psi_2$, the supremum in \eqref{bump_02}, but not on the lattice $\cL$.
	\item For all Haar shifts $\sha_n$ the operators $M_{v^{1/2}} \sha_n M_{w^{1/2}}$  uniformly bounded in $L^2$ by $Cn$, where $C$ is the constant from \textup{\cond1}.
	\item Let $\Pi =\Pi_b$ be a paraproduct such that
\begin{align}
\label{L^infty_carl}
	|J|^{-1} \sum_{I\in\cL: I\subset J} \| \Delta\ci I b\|_\infty^2 |I| \le 1 \qquad \forall J \in\cL.
\end{align}
Then  the operator $M_{v^{1/2}} \Pi M_{w^{1/2}}$ is bounded in $L^2$ by $C$, where again $C$ depends on $\Psi_1$, $\Psi_2$, the supremum in \eqref{bump_02}, but not on the lattice $\cL$.
\end{enumerate}
\end{thm}

\begin{rem}
For the homogeneous lattices, i.e.~ for lattices satisfying
\[
\inf_{J\in\cL} \ \inf_{I\in \ch(J)} \frac{|I |}{|J|} =:\delta >0
\]
all the normalized $L^p$ norms $|I|^{-1/p} \|\Delta\ci I g\|_p$, $p\in [1, \infty]$ are equivalent in the sense of two sided estimates. So for such lattices condition \eqref{L^infty_carl} means that $\|\Pi\|\le C(\delta)$. So Theorem \ref{t:shift-para} gives the estimates that  being fed to Theorem \ref{deco}   imply Theorem \ref{t:main_01}.
\end{rem}
As it was discussed above, \cond 1 implies \cond2. Statement \cond1 is obtained from the following embedding theorem:

\begin{thm}
\label{t:fd-embed}
Let $\Psi$ be as above. Then for any weight $w$ on $X$ such that $\bn\ci \Psi(N_I^w)<\infty$ for all $I\in \cL$
\begin{align}
\label{d-embed-01}
\sum_{I\in\cL} \bn\ci \Psi(N_I^w)^{-1} \biggl(|I|^{-1}\int_X |\Delta\ci I (fw^{1/2})|dx \biggr)^2 |I| \le C\|f\|_{L^2(dx)}^2
\end{align}
for all $f\in L^2(dx)$; here $C=C(\Psi)$ and in the summation we skip $I$ on which $w\equiv0$.

\end{thm}

Let us see that this theorem implies the condition \cond1 of Theorem \ref{t:shift-para}. Assume, multiplying the weights by appropriate constants that the inequality
\begin{align}
\label{bump_03}
\bn\ci{\Psi_1}(N_I^w) \bn\ci{\Psi_2}(N_I^v) \le 1
\end{align}
 holds for all $I\in\cL$. Then
\begin{align*}
|\La \sha (fw^{1/2}), gv^{1/2} \Ra| & \le \sum_{I\in\cL} |\La \sha\ci I \Delta\ci I (fw^{1/2}), \Delta\ci I(gv^{1/2}) \Ra|
\\
& \le
\sum_{I\in\cL} |I|^{-1} \|\Delta\ci I (fw^{1/2}) \|_1 \|\Delta\ci I(gv^{1/2})\|_1
\\
& \le
\sum_{I\in\cL} |I|^{-1} \frac{\|\Delta\ci I (fw^{1/2}) \|_1 \|\Delta\ci I(gv^{1/2})\|_1}{\left( \bn\ci{\Psi_1}(N_I^w) \bn\ci{\Psi_2}(N_I^v)  \right)^{1/2}}
\\ & \le
\frac12 \sum_{I\in\cL} |I|^{-1}\frac{\|\Delta\ci I(fw^{1/2})\|_1^2}{  \bn\ci{\Psi_1}(N_I^w)  }   + \frac12 \sum_{I\in\cL}
|I|^{-1} \frac{  \|\Delta\ci I(gv^{1/2})\|_1^2  }{ \bn\ci{\Psi_2}(N_I^v) }  \,.
\end{align*}
The second inequality here follows from ``$L^1\times L^1$ normalization'' condition \eqref{L^1xL^1_normalization}, the second one from \eqref{bump_03}  and the last one is just the trivial inequality $2xy\le x^2 +y^2$.

Applying Theorem \ref{t:fd-embed} to each sum we get that
\[
|\La \sha (fw^{1/2}), gv^{1/2} \Ra| \le \frac12 \left(  C(\Psi_1) \|f\|_2^2 + C(\Psi_2) \|g\|_2^2    \right)  .
\]
Replacing $f \mapsto tf$, $g\mapsto t^{-1} g$, $t>0$ we get
\[
|\La \sha (fw^{1/2}), gv^{1/2} \Ra| \le \frac12 \left( t^2 C(\Psi_1) \|f\|_2^2 + t^{-2} C(\Psi_2) \|g\|_2^2    \right).
\]
Taking infimum over all $t>0$ and recalling that $2 ab = \inf_{t>0} (t^2a +t^{-2} b)$ for $a, b\ge 0$ we obtain
\[
|\La \sha (fw^{1/2}), gv^{1/2} \Ra| \le (C(\Psi_1)C(\Psi_2))^{1/2} \|f\|_2\|g\|_2,
\]
which is exactly statement \cond1 of Theorem \ref{t:shift-para}. \hfill \qed

For the statement \cond3 of Theorem \ref{t:shift-para} we also need another embedding theorem.
\begin{thm}
\label{t.embed-bump}
Let $\Psi$ be as above.
Then for any normalized Carleson sequence $\{a\ci I\}\ci{I\in\cD}$ ($a\ci I\ge 0$), i.e.~for any sequence satisfying
\[
\sup_{I\in\cD} |I|^{-1}\sum_{I'\in\cD: I'\subset I} a\ci{ I'} |I'| \le 1
\]
we get
\[
\sum_{I\in\cD} \frac{\La fw^{1/2}\Ra\ci I^2}{\bn\ci\Psi(N_I^w)} a\ci I |I| \le C \|f\|_{L^2(dx)}^2,
\]
where again $C=C(\Psi)$.
\end{thm}

Let us show that this theorem together with Theorem \ref{t:fd-embed} implies statement \cond3 of Theorem \ref{t:shift-para}.
Let $a\ci I = \|\Delta\ci I b\|_\infty^2$.

Again, multiplying if necessary the weights $v$ and $w$ by appropriate constants we can assume \eqref{bump_03}.
Then we can write
\begin{align*}
|\La \Pi_b (fw^{1/2}), g v^{1/2} \Ra |
& \le \sum_{I\in\cD} |\La fw^{1/2} \Ra\ci I | \cdot |\La \Delta\ci I b, \Delta\ci I( gv^{1/2}) \Ra |
\\
& \le  \sum_{I\in\cD} \frac{|\La fw^{1/2} \Ra\ci I  | (a\ci I)^{1/2}|I|^{1/2}}{\left( \bn\ci{\Psi_1}(N_I^w) \right)^{1/2} } \cdot
\frac{\| \Delta\ci I( gv^{1/2})\|_1}{\left( \bn\ci{\Psi_2}(N_I^v) \right)^{1/2} |I|^{1/2}}
\\
& \le
\left(
\sum_{I\in\cD} \frac{|\La fw^{1/2} \Ra\ci I  |^2 a\ci I }{\bn\ci{\Psi_1}(N_I^w)  } |I|
\right)^{1/2}
\left(
\sum_{I\in\cD} \frac{\| \Delta\ci I( gv^{1/2})\|_1^2}{\bn\ci{\Psi_2}(N_I^v) |I| }  ;
\right)^{1/2}
\end{align*}
the second inequality holds because of \eqref{bump_03}, and the last one is just the Cauchy--Schwarz inequality.

Estimating the sums in parentheses by Theorem \ref{t.embed-bump} and \ref{t:fd-embed} respectively we get
statement \cond3 of Theorem \ref{t:shift-para}.
\hfill\qed

\section{Proof of (the Differential Embedding) Theorem \ref{t:fd-embed} }

\subsection{Bellman function and main differential inequality}
\label{s:Bell01-DE-01}
Let $\f(s):={s\Psi(s)}$. Multiplying $\Psi$ by an appropriate constant we can assume without loss of generality that
\begin{equation}
\label{norm-Psi}
\int_0^1 \frac1{\f(s)} ds =1.
\end{equation}

 Define $m(s)$ on $[0,1]$ by $m(0)=m'(0)=0$, $m''(s)=1/\f(s)$. Identity \eqref{norm-Psi} implies that $m$ is well-defined and $0\le m'(s) \le 1$,  $0\le m(s)\le s$.
For a distribution function $N=N_I^w$ define
\begin{align}
\label{u(N)}
\bu (N) = \int_0^\infty (2 N(t) - m(N(t))) dt = 2 \La w\Ra\ci I - \int_0^\infty m(N(t)) dt ;
\end{align}
Note that the inequality $m(s) \le s $ implies that $\bu(N_I^w) \ge \La w\Ra\ci I$.

The functional $\bu$ is defined on the convex set of distribution functions, i.e. on the set of decreasing  functions $N:[0, \infty)\to [0,1]$ such that $\int_0^\infty N(t) dt<\infty$.

In what follows we can consider only finitely supported functions $N$, and then use standard approximation reasoning.
Consider two distribution functions $N$ and $N_1$ and let $\sd N= N_1-N$.
Denote also
\[
\bw :=\int_0^\infty N(t) dt, \qquad \bw_1 :=\int_0^\infty N_1(t) dt,
\]
and let
\[
\sd \bw := \bw_1 - \bw = \int_0^\infty \sd N(t) dt;
\]
the motivation for this notation is that if $N$ and $N_1$ are the distribution functions of the weights $w$ and $w_1$, then the integrals are the averages on the corresponding weights.
Denote also
\begin{align}
\label{delta-w}
\bw_\ssd:= \int_0^\infty |\sd N(t)| dt;
\end{align}
clearly $|\sd\bw|\le \bw_\ssd$.

Let us compute derivatives of $\bu$ in the direction of $\sd N$.
The first derivative is given by
\begin{align*}
\bu'\ci{\sd N} (N) =\frac{d}{d\tau}  \bu(N+\tau \sd N)\Bigm|_{\tau=0}
 = \int_0^\infty \left(2  - m'(N(t))\right) \sd N(t) dt ,
\end{align*}
so, in particular
\begin{align*}
|\bu'\ci{\sd N} | \le 2  \bw_\ssd .
\end{align*}
Therefore we can write
\begin{align}
\label{u'}
\bu'\ci{\sd N} = \kappa \bw_\ssd, \qquad \kappa=\kappa(\sd N), \,\,|\kappa|\le 2.
\end{align}

The second derivative in the direction $\sd N =N_1-N$ is given by
\[
-\bu''\ci{\sd N} (N) = -\frac{d^2}{d\tau^2} \bu(N+\tau \sd N)\Bigm|_{\tau=0} =  \int_0^\infty \f(N(t))^{-1}(\sd N(t))^2\,dt
\]
By Cauchy-Schwarz, the integral in the right side is at least
\begin{align*}
\Bigl[\int_0^\infty N(t)\Psi(N(t))\,dt \Bigr]^{-1}
\Bigl[\int_0^\infty |\sd N(t)|\,dt \Bigr]^2\
& = \bn(N)^{-1} \Bigl[\int_0^\infty |\sd N(t)|\,dt \Bigr]^2
\\
& = \bn(N)^{-1} ( \bw_\ssd)^2,
\end{align*}
so
\begin{align}
\label{u''}
-\bu''\ci{\sd N} (N) \ge \frac{(\bw_\ssd)^2}{\bn(N)}
\end{align}

For the scalar variable $f\in \R$ and  the distribution function $N$ define the Bellman function $\wt\cB(f, N)=\cB(\bff, \bu(N))$ where
\[
\cB(\bff,\bu ) = \frac{\bff^2}{\bu}.
\]
Computing second derivative of $\wt\cB$ in the direction $\sd = (\sd \bff, \sd N)$ we get
\begin{align*}
\wt\cB_{\sd}'' =
\left( \begin{array}{c} \sd\bff \\ \bu'_{\sd N} \end{array}\right)^T
\left( \begin{array}{cc} \cB_{\bff \bff} & \cB_{\bff \bu} \\ \cB_{\bff \bu} & \cB_{\bu \bu} \end{array}\right)
\left( \begin{array}{c} \sd\bff \\ \bu'_{\sd N} \end{array}\right) + \cB_\bu \bu''_{\sd N}
\end{align*}

In the last formula the derivative of $\wt\cB$ is evaluated at the point $(f, N)$, and derivatives of $\cB$ are evaluated at $(\bff, \bu(N))$.

The Hessian is easy to compute
\begin{align}
\label{Hess_B-1}
\left( \begin{array}{cc} \cB_{\bff \bff} & \cB_{\bff \bu} \\ \cB_{\bff \bu} & \cB_{\bu \bu} \end{array}\right) =
\left( \begin{array}{cc} \frac{2 }{\bu} & -\frac{2\bff}{\bu^2} \\ -\frac{2\bff}{\bu^2} & \frac{2\bff^2}{\bu^3} \end{array}\right);
\end{align}
note that this matrix is positive semidefinite.

Since $\cB_\bu = -\bff^2/\bu^2$, we get using \eqref{u''}
\[
\cB_\bu \bu''_{\sd N} \ge \frac{\bff^2}{\bu^2\bn}   (\bw_\ssd)^2.
\]
Thus, gathering everything and using  \eqref{u'}  we get
\begin{align}
\label{MainIneq1}
\wt\cB_{\sd}'' \ge
\left( \begin{array}{c} \sd\bff \\ \kappa\bw_\ssd \end{array}\right)^T
\left( \begin{array}{cc} \frac{2 }{\bu} & -\frac{2\bff}{\bu^2} \\ -\frac{2\bff}{\bu^2} & \frac{2\bff^2}{\bu^3} (1+\frac{\bu}{2\kappa^2\bn}) \end{array}\right)
\left( \begin{array}{c} \sd\bff \\ \kappa\bw_\ssd \end{array}\right)
\end{align}
The matrix here is obtained from the Hessian in \eqref{Hess_B-1} by multiplying the lower right entry by $1+\frac{\bu}{2\kappa^2\bn} \ge 1$, so it has more positivity than the Hessian. In particular, if we divide the upper left entry of the matrix in \eqref{MainIneq1} by the same quantity $1+\frac{\bu}{2\kappa^2\bn}$, the matrix still be positive semidefinite. But our matrix in \eqref{MainIneq1} has something bigger in the upper-left corner!

Therefore, since
\[
1- \left( 1 + \frac{\bu}{2\kappa^2\bn}\right)^{-1} = \frac{\bu}{2\kappa^2\bn + \bu}
\]
we get that
\begin{align}
\label{MainIneq2}
\wt\cB_{\sd}'' \ge \frac{2(\sd\bff)^2}{2\kappa^2\bn + \bu} \ge \frac{2(\sd\bff)^2}{2\cdot 2^2\bn + \bu} \ge c \frac{(\sd\bff)^2}{\bn};
\end{align}
the last inequality holds for some $c>0$ because $\bu\le 2 \bw \le C \bn$.

Let us explain it. In fact, we want
$$
\int N_I(t)dt = \av{w}{I} \le C \int N_I(t)\Psi(N_I(t))dt.
$$
Clearly, it is enough to consider the set $B=\{t\colon \Psi(N_I(t)) \le 1\}$. Since $\Psi$ is decreasing, for $t\in B$ we get that $N_I(t)\ge \Psi^{-1}(1)$. Since $s\mapsto s\Psi(s)$ is increasing, we get $N_I(t)\Psi(N_I(t))\ge \Psi^{-1}(1) \ge \Psi^{-1}(1)N_I(t)$ (the last is because $N_I$ is normalized). We are done.

Inequality \eqref{MainIneq2} is exactly what we will use to obtain the Main inequality in difference form in the next section.

\subsection{Main inequality
in the finite difference form}

\subsubsection{Dyadic case}

\begin{lm}
\label{l:MainIneq}
Let
\begin{align*}
\bff = \frac{\bff_1 + \bff_2}{2}, \qquad N(t) = \frac{N_1(t) + N_2(t)}{2}.
\end{align*}
Then
\begin{align}
\label{d-MainIneq}
\frac12 \Bigl( \cB(\bff_1, \bu(N_1))   + \cB(\bff_2, \bu(N_2)) \Bigr) - \cB(\bff, \bu(N)) \ge \frac{c}4   \cdot  \frac{ (\bff_1 -\bff)^2}{\bn(N)}.
\end{align}
where $c$ is the constant from \eqref{MainIneq2}. (Note that $\bff_1 -\bff = \bff-\bff_2$, so we can replace $(\bff_1 -\bff)^2$ in the right side by $(\bff_2 -\bff)^2$)
\end{lm}
\begin{proof}
Notice that
\begin{align}
\label{conc1}
\frac{s_1+s_2}{2} \Psi\left( \frac{s_1+s_2}{2} \right) \ge \frac{s_1+s_2}{2} \Psi\left( s_1+s_2 \right)
\ge\frac12 s_1 \Psi(s_1);
\end{align}
here  the first inequality holds because $\Psi$ is decreasing and the second one because $s\Psi(s)$ is increasing. Of course, we can interchange $s_1$ and $s_2$ in the above inequality.

Let $\sd\bff := \bff_1 -\bff$, $\sd N :=N_1-N$.
Define
\begin{align*}
F(\tau) = \cB(\bff + \tau\sd\bff, \bu(N +\tau\sd N)) + \cB(\bff - \tau\sd\bff, \bu(N -\tau\sd N))
\end{align*}
Taylor's formula together with the estimate  \eqref{MainIneq2} imply that
\begin{align}
\label{d-conv-1}
F(1) - F(0) \ge \frac{c}2 (\sd\bff)^2 \left( \frac1{\bn(N+\tau\sd N)} + \frac1{\bn(N -\tau\sd N)} \right)
\end{align}
for some $\tau\in (0,1)$.

Estimate \eqref{conc1} implies that
\[
\bn(N) \ge \frac12 \bn (N\pm\tau\sd N),
\]
so
\[
\left( \frac1{\bn(N+\tau\sd N)} + \frac1{\bn(N -\tau\sd N)} \right) \ge \frac1{\bn(N)}.
\]
Then it follows from \eqref{d-conv-1} that
\[
F(1) - F(0) \ge \frac{c}2 \cdot \frac{ (\sd\bff)^2}{\bn(N)}.
\]
Recalling the definition of $F$ and dividing this inequality by $2$ we get \eqref{d-MainIneq}.
\end{proof}

\subsubsection{General case}

Let $\f$ and $\wt\cB$ be as above.
\begin{lm}
\label{maindiff2}
Let $\bff, \bff_k \in\R$, $\alpha_k\in \R_+$ and the distribution functions $N$, $N_k$, $k=1, 2, \ldots, n$ satisfy
\[
\bff = \sum_{k=1}^n \alpha_k \bff_k, \qquad N = \sum_{k=1}^n \alpha_k N_k,  \qquad 
\sum_{k=1}^n \alpha_k =1.  \
\]
Then
\begin{align}
\label{eq:maindiff_02}
-\wt\cB(\bff, N) + \sum_{k=1}^n \alpha_k \wt\cB(\bff_k, N_k) \ge \frac{c}{16}\cdot \frac1{\bn(N)} \left( \sum_{k=1}^n \alpha_k | \bff_k -\bff |  \right)^2
\end{align}
\end{lm}

\begin{proof}
The reasoning below is a ``baby version'' of the reasoning used to prove the main estimate (Lemma 6.1) in \cite{TB2011}.

For a weight  $\alpha = \{\alpha_k\}_{k=1}^n$, $\alpha_k\ge0$, let $\ell^p(\alpha)$ be the weighted (finite-dimensional) $\ell^p$ spaces, $\|x\|_{\ell^p(\alpha)}^p =\sum_{k=1}^n \alpha_k |x_k|^p$ ($\ell^\infty(\alpha)$ is just the usual finite-dimensional $\ell^\infty$).

Let $\La\fdot, \fdot \Ra_\alpha$ be the standard duality $\La x, y\Ra_\alpha = \sum_{k=1}^n \alpha_k x_k y_k$.

Define $\be\in \ell^p(\alpha) $, $\be = (1, 1, \ldots, 1)$.

Consider the quotient space $\cX = \ell^1(\alpha)/\spn\{\be\}$. For $x\in \ell^1(\alpha)$ let $x^0:= x- \|\be\|_{\ell^1(\alpha)}^{-1} \La x, \be \Ra_\alpha \be$, so $\sum_{k=1}^n \alpha_k x^0_k =0$. Then
\begin{align}
\label{eq-QuotNorm}
\|x\|\ci{\cX} \le \|x^0\|_{\ell^1(\alpha)} \le 2\|x\|\ci{\cX}.	
\end{align}
Indeed, the first inequality is trivial (follows from the definition of the norm in the quotient space). As for the second one, $|\La x, \be\Ra_\alpha| \le  \|x\|_{\ell^1(\alpha)}$, 
so it follows from the triangle inequality that
\[
\|x^0\|_{\ell^1(\alpha)} \le \|x\|_{\ell^1(\alpha)} + \|\be\|_{\ell^1(\alpha)}^{-1} |\La x, \be \Ra | \cdot \|\be\|_{\ell^1(\alpha)} \le 2 \|x\|_{\ell^1(\alpha)}.
\]
This inequality remains true if one replaces $x$ by $x-\lambda\be$, $\lambda\in\R$, so the second inequality in \eqref{eq-QuotNorm} is proved.

The dual space $\cX^*$ can be identified with s subspace of $\ell^\infty=\ell^\infty(\alpha)$ consisting of $x^*\in \ell^\infty(\alpha)$ such that $\La \be, x^*\Ra_\alpha =0$ (with the usual $\ell^\infty$-norm).

So, for the vector $x= (x_1, x_2, \ldots x_n)$, $x_k =\bff_k -\bff$ (notice that $\La x, \be\Ra_\alpha =0$ there is $\beta =\{\beta_k\}_{k=1}^n$, $|\beta_k |\le 1$ such that $\sum_{k=1}^n \alpha_k \beta_k =0$ and
\[
\sum_{k=1}^n \alpha_k \beta_k (\bff_k-\bff) = \|x\|\ci{\cX} \ge \frac12 \| x\|_{\ell^1(\alpha)} =\frac12  \sum_{k=1}^n \alpha_k |\bff_k -\bff|.
\]

Define $\bff^+$, $\bff^-$, $N^+$, $N^-$ by
\[
\bff^\pm = \sum_{k=1}^n \alpha_k (1 \pm \beta_k) \bff_k, \qquad N^\pm := \sum_{k=1}^n \alpha_k (1 \pm \beta_k) N_k.
\]

By Lemma \ref{l:MainIneq}
\begin{align}
\label{MainIneq_01}
\frac12 \Bigl( \wt\cB(\bff^+, N^+))   +  \wt\cB(\bff^-, N^-)) \Bigr) - \wt\cB(\bff, N) \ge \frac{c}4   \cdot  \frac{ (\bff^+ -\bff)^2}{\bn(N)}
\end{align}
We know that
\begin{align*}
\bff^+ -\bff = \sum_{k=1}^n \alpha_k \beta_k \bff_k = \sum_{k=1}^n \alpha_k \beta_k (\bff_k -\bff) \ge \frac 12
\sum_{k=1}^n \alpha_k |\bff_k -\bff |
\end{align*}
(the second equality holds because $\sum_{k=1}^n \alpha_k\beta_k=0$), so the right side of \eqref{MainIneq_01} is estimated below by
\[
\frac{c}{16} \cdot \frac{1}{\bn(N)} \left( \sum_{k=1}^n \alpha_k | \bff_k -\bff |  \right)^2
\]
Since the function $\wt\cB$ is convex
\begin{align*}
\wt\cB(\bff^+, N^+) &\le  \sum_{k=1}^n \alpha_k (1 + \beta_k) \wt\cB(\bff_k, N_k),  \\
\wt\cB(\bff^-, N^-) &\le  \sum_{k=1}^n \alpha_k (1 - \beta_k) \wt\cB(\bff_k, N_k)
\end{align*}
and adding these inequalities we can estimate above the left side of \eqref{MainIneq_01} by
\[
-\wt\cB(\bff, N) + \sum_{k=1}^n \alpha_k \wt\cB(\bff_k, N_k)  .
\]
\end{proof}

\subsection{From main inequality \eqref{eq:maindiff_02} to Theorem \ref{t:fd-embed}. }
\label{shainik1}
Fix an interval $I^0$ and let $I_k$ be its children. Applying Lemma \ref{maindiff2} with $\bff_k = \La fw^{1/2}\Ra\ci{I_k}$, $N_k = N_{I_k}^w$ and $\alpha_k = |I_k|/|I^0|$ we get denoting $\tilde f := f w^{1/2}$
\begin{align*}
\frac{c}{16}\cdot  \frac{\|\Delta\ci{I^0} \tilde f\|_1^2}{ \bn(N_{I^0}^w) \cdot |I^0| } \le
- |I^0| \wt\cB (\La \tilde f\Ra\ci{I^0}, N_{I^0}^w) + \sum_{I\in\ch(I^0)} |I| \cdot \wt\cB(\La \tilde f \Ra\ci I, N_I^w )
\end{align*}
Applying this formula to all children of $I^0$, then to their children and adding up the inequalities we get after going $n$ generations down  that
\begin{align*}
\frac{c}{16} \sum_{\substack{I\in \ch_k(I^0) \\ 0\le k < n }}
 \frac{\|\Delta\ci{I} \tilde f \|_1^2}{ \bn(N_{I}^w) \cdot |I| }
 &\le
 - |I^0| \wt\cB (\La \tilde f\Ra\ci{I^0}, N_{I^0}^w) + \sum_{I\in\ch_n(I^0)} |I| \cdot \wt\cB(\La \tilde f\Ra\ci I, N_I^w )
 \\
 & \le
 \sum_{I\in\ch_n(I^0)} |I| \cdot \wt\cB(\La \tilde f\Ra\ci I, N_I^w ).
\end{align*}
We know that $\wt\cB(\bff, N) \le C\frac{\bff^2}{\bu(N)}$, and since (see \eqref{u(N)}) $\bu(N_I^w) \ge \La w\Ra\ci I$ we conclude using the Cauchy--Schwarz estimate $|\La fw^{1/2}\Ra\ci I|^2 \le \La |f|^2\Ra\ci I \La w \Ra\ci I$ that
\[
|I| \cdot \wt\cB(\La \tilde f\Ra\ci I, N_I^w ) \le C |I| \frac{ \La fw^{1/2}\Ra\ci I^2}{\La w \Ra\ci I} = C |I| \La |f|^2 \Ra\ci I = C\int_I |f|^2 d\mu.
\]
Therefore, estimating the right side we get
\[
\frac{c}{16} \sum_{\substack{I\in \ch_k(I^0) \\ 0\le k < n }}
 \frac{\|\Delta\ci{I} \tilde f \|_1^2}{ \bn(N_{I}^w) \cdot |I| }
 \le C \int_{I^0} |f|^2 d\mu.
\]
Since the right side does not depend on $n$ we can make $n \to \infty$, and have the sum in the left side over all $I\in\cL$, $I\subset I^0$.

Taking the sum over all $I^0\in \cL_{-m}$ and letting $m\to\infty$ we get conclusion of the theorem. \hfill\qed

\section{Proof of (the Embedding) Theorem \ref{t.embed-bump}. }

\subsection{An auxiliary function}
Let $\Psi$ be the function from Theorem \ref{t.embed-bump}. Define $\f(s) :=s\Psi(s)$.

For the numbers $A\in[1, 2]$, $N\in \R_+$ define
\[
T(A, N) := N \int_0^{N/A} \frac{1}{\f(s)} ds
\]

\begin{lm}
\label{l:prop-T}
The function $T$ is convex and satisfies the differential inequality
\[
-\frac{\partial T}{\p A} \ge \frac14\cdot  \frac{N^2}{\f(N)} .
\]
\end{lm}

\begin{proof}
Differentiating the integral we get
\begin{align}
\label{dT/dA}
-\frac{\partial T}{\p A} = \frac{N^2}{A^2\f(N/A)} \ge \frac14\cdot \frac{N^2}{\f(N)}.
\end{align}
since $\f$ is increasing and $1\le A\le 2$.

To prove the convexity notice that $T$ is linear on the lines $N= k A$, so the Hessian $d^2 T$ degenerates.

Differentiating \eqref{dT/dA} we get
\begin{align*}
\frac{\partial^2 T}{\p A^2} = N^2 \frac{2A \f(N/A) - N\f'(N/A)}{(A^2 \f(N/A))^2}
\end{align*}
Note that the right side is positive if $s\f'(s) < 2 \f(s)$ (because $\f(s)>0$) .

But for our function even a stronger inequality $s\f'(s) \le \f(s)$ is satisfied!  Indeed, since $\f(s) = s\Psi(s) $ is increasing and $\Psi$ is decreasing, then
\begin{align*}
0\le (s\Psi(s))' = \Psi(s) + s \Psi'(s) \le \Psi(s)
\end{align*}
(the second inequality holds because $\Psi $ is decreasing). Multiplying this inequality by $s$ we get $s\f'(s) \le \f(s)$.

Therefore, since $\f(s)>0$, we conclude that $\frac{\p^2 T}{\p A^2}>0$.

But the Hessian $d^2 T$ is singular, and it is an easy exercise in linear algebra to show that a singular Hermitian $2\times 2$ matrix   with a positive entry on the main diagonal is positive semidefinite.
\end{proof}

\subsection{Bellman function and the main differential inequality. }
Let now $N$ be a distribution function, and let
\[
\bT (A, N) = \int_0^\infty T(A, N(t)) \, dt .
\]

As in Section \ref{s:Bell01-DE-01} assume, multiplying $\Psi$ by an appropriate constant,  that
\[
\int_0^1 \frac{1}{\f(s)} ds =1.
\]
Then $T(A, N(t)) \le N(t)$, so
\[
\bT(A, N) \le \int_0^\infty N(t) dt =:\bw =\bw(N).
\]

For $\bff\in\R$, $M\in[0,1]$ and for a distribution function $N$ define the function
$
\wt\cB(\bff, N, M) : = \cB(\bff, \bu(M,N))
$,
where
\[
\cB(\bff, \bu) = \frac{\bff^2}{\bu}
\]
and
\begin{align*}
\bu =\bu(M,N) & = 2 \int_0^\infty N(t) dt - \bT( M+1, N)
\\
& =: 2 \bw(N) - \bT( M+1,N).
\end{align*}
Note that $2 \bw(N) \ge \bu(M,N) \ge \bw(N)$.

We claim that the function $\wt\cB$ is convex. Indeed, fix a direction $\sd := (\sd\bff, \sd N, \sd M)^T$ and compute the second derivative $\wt\cB''_\ssd$ in this direction
\[
\wt\cB''_\ssd = \frac{d^2}{d\tau^2} \wt \cB (\bff+\tau\sd \bff, N+\tau \sd N, M+\tau\sd M)\Bigm|_{\tau=0}.
\]
We get
\[
\wt\cB''_\ssd =
\left(\begin{array}{c} \sd \bff \\ \bu'_\ssd \end{array} \right)^T
\left( \begin{array}{cc} \cB_{\bff\bff} & \cB_{\bff\bu} \\ \cB_{\bff\bu} & \cB_{\bu\bu} \end{array} \right)
\left(\begin{array}{c} \sd \bff \\ \bu'_\ssd \end{array} \right)
+\cB_{\bu} \bu''_\ssd.
\]
The Hessian
\[
\left( \begin{array}{cc} \cB_{\bff\bff} & \cB_{\bff\bu} \\ \cB_{\bff\bu} & \cB_{\bu\bu} \end{array} \right)
=
\left( \begin{array}{cc} \frac{2 }{\bu} & -\frac{2\bff}{\bu^2} \\ -\frac{2\bff}{\bu^2} & \frac{2\bff^2}{\bu^3} \end{array}
\right)
\]
is clearly positive semidefinite, so the first term is nonnegative. For the second term notice that
\begin{align}
\label{cB_bu}
\cB_{\bu} = -\frac{\bff^2}{\bu^2}, \qquad \bu''_\ssd = -\bT''_\ssd \le 0
\end{align}
because $T$, and therefore $\bT$ is convex. Thus $\wt\cB''_\ssd\ge 0$, so $\wt\cB$ is convex.
Let us compute the partial derivative
\begin{align}
\label{dB/dM-01}
-\frac{\p \wt\cB}{\p M} = -\cB_\bu \frac{\p\bu}{\p M} = \frac{\bff^2}{\bu^2}\cdot \left( - \frac{\p \bT}{\p M}\right)
\end{align}
By Lemma \ref{l:prop-T}
\begin{align*}
- \frac{\p \bT}{\p M}  & \ge \frac14 \cdot \int_0^\infty \frac{N(t)^2}{\f(N(t))} dt
\\
& \ge \frac14 \left(\int_0^\infty N(t) dt \right)^2 \left( \int_0^\infty \f(N(t)) dt \right)^{-1}
= \frac14 \cdot \frac{\bw(N)^2}{\bn(N)};
\end{align*}
the second inequality here is just the Cauchy--Schwarz inequality. Combining with \eqref{dB/dM-01} and recalling that $\bu\le 2\bw$ we get
\begin{align}
\label{dB/dM-02}
-\frac{\p \wt\cB}{\p M} \ge \frac1{16}\cdot \frac{\bff^2}{\bn}
\end{align}
This inequality (together with the convexity of $\wt\cB$) is the main differential inequality for our function.

\subsection{Finite difference form of the main inequality}
Let $X=(\bff, N, M)$, $X_k =(\bff_k, N_k, M_k)$, ($\bff, \bff_k\in \R$, $M, M_k\in[0,1]$, $N$, $N_k$ are the distribution functions) satisfy
\begin{align*}
\bff = \sum_{k=1}^n \alpha_k \bff_k, \qquad N = \sum_{k=1}^n \alpha_k N_k, \qquad M = a + \sum_{k=1}^n \alpha_k M_k, \ a\ge 0,
\end{align*}
where
\begin{align*}
 \sum_{k=1}^n \alpha_k =1, \qquad \alpha_k \ge 0.
\end{align*}
Then
\begin{align}
\label{discr-MainIneq2-03}
- \wt \cB(X) + \sum_{k=1}^n \alpha_k \wt\cB(X_k) \ge \frac1{16} \cdot \frac{a\bff^2}{\bn}
\end{align}
where $\bn=\bn(N)$.

Indeed, for $M_0:= \sum_{k=1}^n \alpha_k M_k$ the main inequality \eqref{dB/dM-02} implies
\begin{align*}
\wt\cB(\bff, N, M_0) - \wt\cB(\bff, N, M) \ge \frac1{16} \cdot \frac{a\bff^2}{\bn}.
\end{align*}
The convexity of $\wt\cB$ implies that
\[
\wt\cB(\bff, N, M_0) \le \sum_{k=1}^n \alpha_k \wt\cB(X_k)
\]
which together with the previous inequality gives us \eqref{discr-MainIneq2-03}.

\subsection{From main inequality \eqref{discr-MainIneq2-03} to Theorem \ref{t.embed-bump}.}

The reasoning here is almost verbatim the same as in Section \ref{shainik1}.

For an interval let us $I\in\cL$ denote $\bff\ci I=\La fw^{1/2}\Ra\ci I$, $N_I =N_{I}^w$, $M\ci{I}= |I|^{-1}\sum_{I'\subset I} a\ci{I'}$, $\bw\ci I =\La w \Ra\ci I$, $\bu\ci I = \bu(M\ci I, N\ci I)$.

Fix  $I^0\in\cL$, and let $I_k$ be its children.
Applying the inequality \eqref{discr-MainIneq2-03} with $\alpha_k = |I_k|/|I^0|$, $\bff_k = \bff\ci{I_k}$, $N_k = N_{I_k}^w$, $M_k= M_{I_k}$ we get that
\[
\frac1{16} \cdot \frac{a\ci{I^0} \bff\ci{I^0}^2}{\bn(N_{I^0}^w)} |I^0| \le
-|I^0| \wt\cB (X\ci{I^0}) + \sum_{I\in\ch(I^0)} |I|\wt\cB(X\ci I)
\]
Writing the corresponding estimates for the children of $I^0$, then for their children, we get after going $n$ generations down and using the telescoping sum in the right side
\begin{align*}
\frac1{16}\sum_{\substack{ I\in \ch_k (I^0)\\ 0\le k <n}} \frac{a\ci{I} \bff\ci{I}^2}{\bn(N_{I}^w)} |I|
& \le
-|I^0| \wt\cB (X\ci{I^0}) + \sum_{I\in\ch_n(I^0)} |I|\wt\cB(X\ci I)
\\ &
\le \sum_{I\in\ch_n(I^0)} |I|\wt\cB(X\ci I);
\end{align*}
the last inequality holds because $\wt\cB\ge 0$.

Since
\[
\wt\cB(X\ci I) \le \bff\ci I^2/\bu\ci I \le \bff\ci I^2/\bw\ci I
\]
(the last inequality holds because $\bu\ge \bw$) and by Cauchy--Schwarz
\[
|\La fw^{1/2}\Ra\ci I|^2 \le \La |f|^2\Ra\ci I \La w \Ra\ci I,
\]
we conclude, exactly as in Section \ref{shainik1} that
\[
|I|\wt\cB(X\ci I) \le |I|\La|f|^2\Ra\ci I = \int_I |f|^2 d\mu,
\]
so
\begin{align*}
\frac1{16}\sum_{\substack{ I\in \ch_k (I^0)\\ 0\le k <n}} \frac{a\ci{I} \bff\ci{I}^2}{\bn(N_{I}^w)} |I|
\le \int_{I^0} |f|^2 d\mu.
\end{align*}
Conclusion of the proof is exactly as in Section \ref{shainik1}: we first let $n\to \infty$, and then taking the sum over $I^0\in \cL_{-m}$ and letting $m\to\infty$ get the desired estimate. \hfill\qed

\section{Concluding remarks and open problems}

\subsection{One sided bumps}
The famous theorem of P.~Koosis \cite{Koosis-1/w-1980}  states that given a weight $u$ on the unit circle, one can find a non-zero weight $v$ such that
the Hilbert transform $T$ is  bounded as an operator from $L^2(u)$ to $L^2(v)$ if and only if $1/u\in L^1$.
The same result hods for the maximal function, see \cite{RDF-KoosisThm_1981}.

So it is possible to have a situation when one has two weight estimates, but one cannot ``bump'' the $L^1$-norm of $1/u$ ($w$ in our notation). This leads to a very natural in our opinion) question ``can one ``bump'' only one weight to get the two weight estimate?'' For example, can one find a reasonable Young function $\Phi$ that the condition
\[
\sup_{I} \left(|I|^{-1} \int_I w \right) \|v\|\ci{L^\Phi(I)} <\infty
\]
implies the boundedness in$L^2$ of the operator $M_{v^{1/2}} T M_{w^{1/2}}$ for all (bounded) \cz operators $T$?

Note, that for maximal function a one sided bump  condition (but with the bump on the ``wrong" side) is sufficient. Namely, it follows from the result in \cite{PeMax} that
if $\int^\infty (\Phi(x))^{-1} dx<\infty$ and the weights $v$, $w$ satisfy
\begin{align}
\label{bump_max1}
\sup_I \left(|I|^{-1} \int_I v \right) \| w\|\ci{L^{\Phi}(I)} < \infty\,,
\end{align}
then the operator $M_{v^{1/2}} \cM M_{w^{1/2}}$, where $\cM$ is the Hardy--Littlewood maximal operator, is bounded in $L^2$. It is natural to remark here that in \cite{NRTV1} we demonstrated this result of P\'erez by almost precisely the same Bellman function that the reader saw above.



\subsection{Estimates for general measures}
A standard way to set up the two weight estimate problem for integral operators is to make the change of variables so  in the integral operator one integrates with respect to the same measure that is used to compute the norm in the domain.

Namely, if one defines measures $\mu$, $\nu$, $d\mu= wdx$, $d\nu=vdx$, then the $L^p$ boundedness of the operator  $M_{v^{1/p} }T M_{w^{1/p'}}$ is equivalent (at least formally) to the boundedness of the operator $T_\mu :L^p(\mu)\to L^p(\nu)$, where
\[
T_\mu f (x) =\int K(x, y) f(y) d\mu(y);
\]
$K$ here is the kernel of the \cz operator $T$.

In fact, everything can be interpreted absolutely rigorously. The boundedness of such operators can be interpreted  as uniform boundedness of the smooth truncations; in  fact such uniform boundedness is equivalent to the boundedness of the bilinear form on functions with separated compact supports, i.e.~to the weakest possible notion of boundedness, see \cite{Liaw-Tr_Reg_2010}.

This setting gives the most general form of the two weight problem, since  $\mu$ and $\nu$ can be general measures, not necessarily absolutely continuous, they even can be purely singular. So the question arises, ``how one can bump general measures?'' The approach with Orlicz spaces, or any other function spaces works only for absolutely continuous measures.

\def\cprime{$'$}
  \def\lfhook#1{\setbox0=\hbox{#1}{\ooalign{\hidewidth\lower1.5ex\hbox{'}\hidewidth\crcr\unhbox0}}}
\providecommand{\bysame}{\leavevmode\hbox to3em{\hrulefill}\thinspace}
\providecommand{\MR}{\relax\ifhmode\unskip\space\fi MR }
\providecommand{\MRhref}[2]{%
  \href{http://www.ams.org/mathscinet-getitem?mr=#1}{#2}
}
\providecommand{\href}[2]{#2}

\end{document}